\documentclass[a4paper]{article}
 \pdfoutput=1          %% arxiv在线编辑需要
 \UseRawInputEncoding  %% arxiv在线编辑需要 该行在本地运行时候必须删除
\usepackage[textwidth=14.5cm]{geometry}  %解决行溢出的问题
\usepackage{blindtext}                   %解决行溢出的问题
\usepackage{amsmath,amsfonts,amssymb,amsthm}
\usepackage{latexsym}
\usepackage{enumerate}
\usepackage{graphics}
\usepackage{epsfig}
\usepackage{mathdots}
\usepackage{caption}
\usepackage{subcaption}
\usepackage{setspace}
\usepackage{multirow}
\usepackage{dcolumn}
\usepackage{url}
\usepackage{tikz,pgfplots}
\usetikzlibrary{arrows,backgrounds}
\usepackage{ragged2e}
\usepackage[numbers,sort&compress]{natbib}

\newtheorem{thm}{Theorem}
\newtheorem{lem}{Lemma}

\newtheorem{defn}{Definition}  %%% zyw add
\newtheorem{exmp}{Example}  %%% zyw add
\newenvironment{solution}{\emph{Solution.}}  %% zyw add

\author{Yongwen Zhu\\
\small School of Mathematics and Information Science, Yantai University,\\Yantai City 264005, P.R. China
\\ \small Email: {zyw@ytu.edu.cn}
}

\title{A Multiplication Formula and Its Application
}

\date{}

\begin{document}

\begin{sloppypar}     %解决行溢出的问题 solve overfull\hbox

\maketitle

\abstract{
This article presents a vertical multiplication formula for calculating the multiplication of any two multi-digit integers, which may be not only used to design the multiplier but also to the mental multiplication.
Our algorithm is a generalization of that of Karatsuba and Offman, but it is  superior to the latter because the latter is not suitable for oral calculation of multiplication of large integers. These two algorithms are both the vertical multiplication method, but in our algorithm the horizontal subtraction in stead of the horizontal addition used in the Kartsuba algorithm makes the involved numbers are much smaller and their algebraic sum are much smaller too because the positive and negative numbers probably cancel each other out. Because of the above advantages, our vertical multiplication formula allows more practical and efficient multipliers to be designed.
}


%\textbf{MR(2010) Subject Classification} 11A99, 20N99%%  11A  Elementary number theory, 11Nxx 乘法数理论 20M半群

\section{Introduction}

Multipliers requiring large bit lengths have a major impact on the performance of many applications, such as cryptography, digital signal processing (DSP) and image processing. As multipliers take a long time for execution so there is a need of fast multiplier to save the execution time. Vedic algorithms are used to the design of computer processors for enhancing speed and performance, see
\cite{Bansal,Garg,Gupta,Gurumurthy,Kavita,Kayal,Mathur,Paramasivam,Pradhan,Ramalatha,Rashno,Sahu}.

 Novel, optimised designs of large integer multiplication are needed as previous approaches, such as schoolbook multiplication, may not be as feasible due to the large parameter sizes. Parameter bit lengths of up to millions of bits are required for use in cryptography, such as in lattice-based and fully homomorphic encryption (FHE) schemes. In \cite{Rafferty},  Rafferty et al presented a comparison of hardware architectures for large integer multiplication, and several multiplication methods and combinations thereof are analyzed for suitability in hardware designs, targeting the FPGA platform.  In particular, the first hardware architecture combining Karatsuba and Comba multiplication was proposed by them.

In \cite{San}, San et al presented a compact hardware architecture for long integer multiplication and proposed a strategy to increase the computational efficiency of the Karatsuba algorithm on FPGA, which provides an efficient and compact architecture to be used where long integer multiplication is definitely required such as Cryptography, especially Public Key Cryptography (PKC), Coding theory, DSP and many more.

Fully homomorphic encryption (FHE) is a technique enabling processing to be performed directly on encrypted data in a commercial cloud environment, thereby preserving privacy. Large integer multiplication is the most time-consuming operation during the FHE.
In \cite{Hua}, Hua et al  proposed an operands merging method of the number theoretic transform (NTT) multiplication butterfly unit. Then by using the operands merging method and a fast modulo method, the multiplier would improve the area requirements and performance. For a combined single trace attack on global shuffling long integer multiplication and its novel countermeasure, one can consult \cite{Lee}.

In \cite{Shi-Wang}, Shi W. et al presented a configurable long integer multiplier tailored to subthreshold operation for ultra-low-power asymmetric en/decryption in semi-passive or passive systems. Their multiplier was composed of radix-4 booth de/encoders and two $544\times 164$-bit partial product reduction tree arrays and reconfigurable data paths, which can be configured for multiple multiplications with different bit size.

Recently, Harvey and Hoeven presents an algorithm that computes the product of two $n$-bit integers in $O(n \log n)$ bit operations, thus confirming a conjecture of Sch{\"o}nhage and Strassen from 1971.  This new algorithm adopts a novel “Gaussian resampling” technique that reduces the integer multiplication problem to a collection of multidimensional discrete Fourier transforms over the complex numbers, which may be evaluated rapidly by means of Nussbaumer' s fast polynomial transforms, See \cite{Harvey-Hoeven}.
However, in its current form, the new algorithm is not actually practical because the proof presented in the paper only works for very large numbers.
Even if every number were written on hydrogen atoms, there would hardly be enough room to write them down in the observable universe.

In the present paper, we present a multiplication formula, which may be used to design the multiplier but also to mental multiplication.
Our algorithm is a generalization of that of Karatsuba and Offman
\cite{Rafferty,Harvey-Hoeven,Kavita,Karatsuba}, but it is  superior to the latter because the latter is not suitable for oral calculation of multiplication. Two algorithm are both the vertical multiplication method, but in our algorithm the horizontal subtraction in stead of the horizontal addition of Kartsuba algorithm makes the involved numbers are much smaller and their algebraic sum are much smaller too because the positive and negative numbers probably cancel each other out.

\section{The basic formula of the rapid multiplication of two integers}

As preparation of our main formula, we present the basic formula of the rapid multiplication of two integers, which was in effect used by some authors, see for example, \cite{Kavita,Shi}. But here is the most specific and most general form of the formula.
\par
Let's begin with a notion and its notation. For any two groups of numbers $a_1,a_2,\ldots,a_n$ and $b_1,b_2,\ldots, b_n$, the \emph{cross product sum} is defined as
$$\begin{bmatrix} a_1 & a_2 & \cdots & a_n \\ b_1 & b_2 & \cdots & b_n  \end{bmatrix}=\sum_{i=1}^n a_i\times b_{n+1-i}=a_1\times b_n+a_2\times b_{n-1}+\cdots+a_n\times b_1.$$
For example, the cross product sum of $1,2,3$ and $4,5,6$ is $$\begin{bmatrix} 1 & 2 &3 \\4 & 5  & 6 \end{bmatrix}=1\times 6 +2\times 5 +3\times 4=28.$$
Evidently, the cross product sum of two numbers is exactly their product, that is,
$$\begin{bmatrix} a\\b\end{bmatrix}=a\times b.$$
\par
For example, to calculate the product $123\times 456$, we can use the cross product sums.
If we momentarily ignore the carry of all involved cross product sums, then the units of this  product is equal to $3\times 6$, that is $\begin{bmatrix} 3\\6 \end{bmatrix}$;
the tens is $2\times 6+3\times 5$, that is $\begin{bmatrix} 2 & 3\\5 & 6 \end{bmatrix}$;
the hundreds is $1\times 6+2\times 5+3\times 4$, that is $\begin{bmatrix} 1 & 2 & 3\\4 & 5 & 6 \end{bmatrix}$;
the thousands is $1\times 5+2\times 4$, that is $\begin{bmatrix} 1 & 2 \\4 & 5 \end{bmatrix}$;
the ten-thousands is $1\times 4$, that is $\begin{bmatrix} 1 \\4\end{bmatrix}$.
Therefore, $$123\times 456=(\begin{bmatrix} 1 \\4\end{bmatrix},\begin{bmatrix} 1 & 2 \\4 & 5 \end{bmatrix},\begin{bmatrix} 1 & 2 & 3\\4 & 5 & 6 \end{bmatrix},\begin{bmatrix} 2 & 3\\5 & 6 \end{bmatrix},\begin{bmatrix} 3\\6 \end{bmatrix}).$$
The final result of the product can be obtained by calculating each cross product sum and carrying out the rounding process:
$$123\times 456=(4,13,28,27,18)=(5,5,t,8,8)=56\,088,$$
where $t$ represents the number $10$, which causes a further carry.
\par

Generally, we have
\begin{lem}[The Basic Formula of the Rapid Multiplication ]\label{lem:chengjijibengongshi-1}\label{basic-formula}
For $m\geq n$, the product of an $m$-digit number $(a_1,a_2,\ldots, a_m)$ and an $n$-digit number $(b_1,b_2,\ldots,b_n)$ equals to
\begin{gather*}
(a_1,a_2,\ldots, a_m)\times (b_1,b_2,\ldots,b_n)\\
=(\begin{bmatrix} a_1 \\b_1\end{bmatrix},\begin{bmatrix} a_1 & a_2 \\b_1 & b_2 \end{bmatrix},\cdots, \begin{bmatrix} a_1 & a_2 & \ldots & a_n\\b_1 & b_2 & \ldots & b_n\end{bmatrix},
\begin{bmatrix} a_2 & a_3 & \ldots & a_{n+1}\\b_1 & b_2 & \ldots & b_n\end{bmatrix},
\ldots,
\\
\begin{bmatrix} a_{m-n+1} & a_{m-n+2} & \ldots & a_m\\b_1 & b_2 & \ldots & b_n\end{bmatrix},
\begin{bmatrix} a_{m-n+2} & a_{m-n+3} & \ldots & a_m\\b_2 & b_3 & \ldots & b_n\end{bmatrix},
\cdots,
\begin{bmatrix} a_{m-1} & a_m \\b_{n-1} & b_n \end{bmatrix},\begin{bmatrix} a_m\\b_n \end{bmatrix}).
\end{gather*}
\end{lem}

For multi-digit numbers, we may use the segmented numbers. For example, for the six-digit number $123\,456$, the segmentations by length $2$ are $(12,34,56)$, and the segmentations by length $3$ are $(123,456)$. Notice that, the above basic formula (in Lemma~\ref{lem:chengjijibengongshi-1}) is also applicable for segmented numbers whose segment length is greater than $1$. Here is an example to use the basic formula of the rapid multiplication to segmented numbers:

\begin{exmp}\rm
Using Lemma~\ref{lem:chengjijibengongshi-1} to segmented numbers with segment length $2$, we can compute the following multiplication in mind:
\begin{gather*}
2\,976\times 2\,924=(29,76)\times (29,24)=
(\begin{bmatrix} 29  \\ 29 \end{bmatrix},\begin{bmatrix} 29 & 76 \\ 29 & 24 \end{bmatrix},\begin{bmatrix}  76 \\ 24 \end{bmatrix})\\
=(\begin{bmatrix} 29  \\ 29 \end{bmatrix},\begin{bmatrix} 29 \\ 76+24 \end{bmatrix},\begin{bmatrix}  76 \\ 24 \end{bmatrix})
=(\begin{bmatrix} 29  \\ 29 \end{bmatrix},\begin{bmatrix} 29 & 100 \end{bmatrix},\begin{bmatrix}  76 \\ 24 \end{bmatrix})
=(\begin{bmatrix} 29  \\ 29+1 \end{bmatrix},0,\begin{bmatrix}  76 \\ 24 \end{bmatrix})\\
=(\begin{bmatrix} 29  \\ 30 \end{bmatrix},0,\begin{bmatrix}  76 \\ 24 \end{bmatrix})
=(\begin{bmatrix} 29  \\ 30 \end{bmatrix},0,\begin{bmatrix}  76 \\ 24 \end{bmatrix})
=(870,00,1824)=(8,70,18,24)=8\,701\,824.
\end{gather*}
\end{exmp}

\section{The vertical multiplication formula}

As the main result of this paper, we shall present and prove the vertical multiplication formula, which may be applied to fast multiplication as well as computer science.

For each cross product sum  $\begin{bmatrix} a_1 & a_2 \\b_1 & b_2 \end{bmatrix}$ with length $2$, the \emph{difference product} is defined as
\[
K_{12} =
\begin{bmatrix} a_1 - a_2 \\b_1 - b_2 \end{bmatrix}
=(a_1-a_2 )\times (b_1 -b_2 )=(a_2-a_1 )\times (b_2 -b_1 ).
\]
Then for any cross product sum
\[\begin{bmatrix} a_1 & a_2 & \ldots & a_n\\b_1 & b_2 & \ldots & b_n\end{bmatrix}\]
with length $n$,
the \emph{symmetric difference} is defined as
\[ K_{12\cdots n} =
\begin{cases}
K_{1n}+K_{2,n-1}+\cdots+K_{m,m+1}, & if n=2m；\\
K_{1n}+K_{2,n-1}+\cdots+K_{m,m+2}, & if  n=2m+1.
\end{cases} \]
Since
\[
K_{12} =\begin{bmatrix} a_1 - a_2 \\b_1 - b_2 \end{bmatrix}
=\begin{bmatrix} a_1  \\b_1 \end{bmatrix} +\begin{bmatrix} a_2  \\ b_2 \end{bmatrix}-\begin{bmatrix} a_1 & a_2 \\b_1 & b_2 \end{bmatrix},
\]
we have that
\begin{gather*}
 K_{12\cdots n}=K_{1n}+K_{2,n-1}+\cdots  \\
=\begin{bmatrix} a_1  \\b_1 \end{bmatrix} +\begin{bmatrix} a_n  \\ b_n \end{bmatrix}-\begin{bmatrix} a_1 & a_n \\b_1 & b_n \end{bmatrix}
+\begin{bmatrix} a_2  \\b_2 \end{bmatrix} +\begin{bmatrix} a_{n-1}  \\ b_{n-1} \end{bmatrix}-\begin{bmatrix} a_2 & a_{n-1} \\b_2 & b_{n-1} \end{bmatrix}
+\cdots  \\
=\begin{bmatrix} a_1  \\b_1 \end{bmatrix}+\begin{bmatrix} a_2  \\b_2 \end{bmatrix} +\cdots+\begin{bmatrix} a_n  \\ b_n \end{bmatrix}-  \begin{bmatrix} a_1 & a_2 & \ldots & a_n\\b_1 & b_2 & \ldots & b_n\end{bmatrix}.
\end{gather*}
Therefore, we obtain the following proposition.

\begin{lem}
Each cross product sum equals to the corresponding vertical product sum minus the symmetric difference product sum, that is,
\[
\begin{bmatrix} a_1 & a_2 & \ldots & a_n\\b_1 & b_2 & \ldots & b_n\end{bmatrix}
=C_1+C_2+\cdots+C_n-K_{12\cdots n},
\]
where $C_i=\begin{bmatrix} a_i \\b_i\end{bmatrix}$ is the vertical product, for $i=1,2,\ldots, n$.
\end{lem}

Let $C=(C_1,C_2,\ldots,C_n)$. $C$ is also called the vertical product of $a$ and $b$.
For any two $n$-digit natural numbers, we apply the the basic rapid multiplication formula~\ref{basic-formula} and the above lemma to obtain that
\begin{gather*}
(a_1,a_2,\ldots, a_n)\times (b_1,b_2,\ldots,b_n)\\
=(\begin{bmatrix} a_1 \\b_1\end{bmatrix},\begin{bmatrix} a_1 & a_2 \\b_1 & b_2 \end{bmatrix},\ldots, \begin{bmatrix} a_1 & a_2 & \ldots & a_n\\b_1 & b_2 & \ldots & b_n\end{bmatrix},
\begin{bmatrix} a_2 & a_3 & \ldots & a_n\\b_2 & b_3 & \ldots & b_n\end{bmatrix},
\cdots,
\begin{bmatrix} a_{n-1} & a_n \\b_{n-1} & b_n \end{bmatrix},\begin{bmatrix} a_n\\b_n \end{bmatrix})\\
=(C_1, C_1+C_2 -K_{12}, C_1+C_2 +C_3-K_{123},\ldots, C_1+C_2 +\cdots+C_n -K_{12\cdots n}, \\
C_2+C_3+\cdots+C_n -K_{2\cdots n},\ldots, C_{n-1}+C_n -K_{n-1,n},C_n)  \\
=(C_1, C_1+C_2, C_1+C_2+C_3,\ldots, C_1+C_2 +\cdots+C_n, C_2+C_3+\cdots+C_n,\ldots, \\
C_{n-1}+C_n,C_n)
-(0, K_{12}, K_{123}, \ldots, K_{12\cdots n}, K_{2\cdots n},\ldots, K_{n-1,n},0)\\
=(C_1,C_2,\ldots,C_n)\times (\underbrace{1,1,\ldots,1}_n)-K,
\end{gather*}
where $K=(0, K_{12}, K_{123}, \cdots, K_{12\cdots n}, K_{2\cdots n},\ldots, K_{n-1,n},0)$ is called the \emph{tare}. Therefore, we obtain the main results of this paper

\begin{thm}[Vertical Multiplication Formula]
Suppose that $K$ is the above tare. Then
\[
(a_1,a_2,\ldots, a_n)\times (b_1,b_2,\ldots,b_n)=(C_1,C_2,\ldots,C_n)\times (\underbrace{1,1,\ldots,1}_n)-K.
\]
\end{thm}

\noindent {\bf Note}： This formula is also able to be used for segmented numbers with each segment length greater than $1$.

\section{A brief introduction of the scissor product theory for multiplication of integers}

Based on the basic formula of the rapid multiplication, in \cite{Zhu}, we proposed the the scissor product theory for mental multiplication as a generalization of that of Shi \cite {Shi}.

\begin{defn}\rm \cite{Zhu}
The difference of a positive integer $a$ and the standard number $N = 10 ^n$ is called the \emph{increment} of $a$ (relative to $N$), denoted as $a^\Delta$, namely, $a^\Delta = a - N$; and the difference of the standard number $N$ and the number $a$ is called the \emph{complement} of $a$ (relative to $N$), denoted as $\widetilde{a}$, that is, $\widetilde{a}=N-a$.
\end{defn}

Obviously, the increment and the complement of $a$ relative to $N$ are opposite numbers reciprocally, i.e.,
$$\widetilde{a}+a^\Delta=0.$$
For example, the increment of $9$ relative to $10$ is $-1$, while the complement of  $9$ relative to $10$ is $1$;
the increment of $119$ relative to $100$ is $19$, while the complement of  $119$ relative to $100$ is $-19$.

\begin{defn}\rm \cite{Zhu}  \label{def:jiandaoji}
Let  $a$, $b$ be two positive integers and $N=10^n$ the standard number. Then the scissor product $a\wedge b$ of $a$ and $b$ relative to $N$ is defined as：
$$a\wedge b=a\times b-(\min(a,b)-1)\times N,$$
where $\min(a,b)$ represents the smaller one of two numbers $a,\,b$. If $a\leq b$, then $a\wedge b$ is also denoted as $\begin{bmatrix} \widehat{a} \\ b \end{bmatrix}$;
if $a\geq  b$, then $a\wedge b$ is also written as $\begin{bmatrix} a \\ \widehat{b} \end{bmatrix}$.
\end{defn}

For example, the scissor product of $3$ and $9$ relative to $10$ is
$$3\wedge 9=\begin{bmatrix} \widehat{3} \\ 9 \end{bmatrix}=3\times 9-(3-1)\times 10=27-20=7,$$
and the scissor product of $2$ and $97$ relative to $100$ is
$$2\wedge 97 =\begin{bmatrix} \widehat{2} \\ 97 \end{bmatrix}=2\times 97-(2-1)\times 100=194-100=94.$$

By the symmetry, we know that the scissor product satisfies the commutative law,  that is,
\begin{lem}\cite{Zhu}
$a\wedge b =b \wedge a. $
\end{lem}

Since $b^\Delta=b-N$,  $N=b-b^\Delta$. So, $a\times N=a\times (b-b^\Delta)=a\times b-a\times b^\Delta$, i.e.,
$a\times b=a\times N+a\times b^\Delta=(a, a\times b^\Delta)$, where the last number is the segmented number with the segment length $n$. Thus, we obtain

\begin{lem}\cite{Zhu}\label{lem:1-zengliang}
$\begin{bmatrix} a \\ b \end{bmatrix}=(a,\begin{bmatrix} a \\ b^\Delta \end{bmatrix}).$
\end{lem}

For the convenience, set $$a^+=a+1,\quad a^-=a-1.$$
According to Lemma \ref{lem:1-zengliang}, we immediately have：
\begin{thm}\cite{Zhu}\label{thm-jiandaoji}
If $a\leq b$, then
$$\begin{bmatrix} a \\ b \end{bmatrix}=(a^-,\begin{bmatrix} \widehat{a} \\ b \end{bmatrix}).$$
\end{thm}

\section{The plum-blossom product method of multiplication of multi-digit numbers}

Note that the absolute value of the scissor product of some digits is not less than $10$, for example, $5\wedge 5=-15$. This would bring us non-convenience when we use scissor product to calculate the multiplication of two multi-digit numbers. For improving the multiplication method of scissor products, here we suggest a new concept called plum-blossom product which is defined below:

\begin{defn}
For two digits $a$ and $b$, the plum-blossom product of $a$ and $b$ is denoted as $a\clubsuit b$
and it is defined as the ones of the usual product of $a$ and $b$ if the ones is less than or equal to $3$ or otherwise the ones minus $10$.
\end{defn}

By the definition, the plum-blossom product of any two digits is less than or equal to $3$ and meanwhile it is greater than or equal to $-6$.
For example, since $3\times 7=21$ and $1\leqslant 3$, the plum-blossom product of $3$ and $7$ is $1$, that is, $3\clubsuit 7=1$. Since  $7\times 7=49$ and $9\geq 3$,  the plum-blossom product of $7$ and $7$ is $9-10=-1$, that is, $7\clubsuit 7=-1$. Similarly, $5\clubsuit 5=5-10=-5$ and $5\clubsuit 8=0$. From the usual multiplication table, one may easily deduce the following plum-product table:

\begin{table}[!htb]
\centering
 \begin{tabular} %%|p{1cm}<{\centering}|
 {|c||c|c|c|c|c|c|c|c|c|c|}
                    \hline
                    $\clubsuit$  &  \rule{8pt}{0pt}$1$\rule{8pt}{0pt}&  \rule{8pt}{0pt}$2$\rule{8pt}{0pt}&  \rule{8pt}{0pt}$3$\rule{8pt}{0pt}&  \rule{8pt}{0pt}$4$\rule{8pt}{0pt}&  \rule{8pt}{0pt}$5$\rule{8pt}{0pt}&  \rule{8pt}{0pt}$6$\rule{8pt}{0pt}&  \rule{8pt}{0pt}$7$\rule{8pt}{0pt}&  \rule{8pt}{0pt}$8$\rule{8pt}{0pt}&  \rule{8pt}{0pt}$9$\rule{8pt}{0pt}\\
                                   \hline\hline
                   $1 $&$ 1 $&$ 2 $&$ 3 $&$ -6 $&$ -5 $&$ -4 $&$ -3 $&$ -2 $&$ -1 $\\
                   \hline
                   $2 $&$  $&$ -6 $&$ -4 $&$ -2 $&$ 0 $&$ 2 $&$ -6 $&$ -4 $&$ -2 $\\
                   \hline
                  $ 3 $&$  $&$  $&$ -1 $&$ 2 $&$ -5$ &$ -2 $&$ 1 $&$ -6 $&$ -3 $\\
                   \hline
                  $ 4 $&$  $&$  $&$  $&$ -4 $&$ 0 $&$ -6 $&$ -2 $&$ 2 $&$ -4 $\\
                   \hline
                   $5 $&$  $&$  $&$  $&$     $&$-5 $&$ 0 $&$ -5 $&$ 0 $&$ -5 $\\
                   \hline
                   $6 $&$  $&$  $&$  $&$  $&$  $&$ -4 $&$ 2 $&$ -2 $&$ -6 $\\
                   \hline
                   $7 $&$  $&$  $&$  $&$  $&$  $&$  $&$ -1  $&$ -4 $&$  3 $\\
                   \hline
                  $ 8 $&$  $&$  $&$  $&$  $&$  $&$  $&$  $&$ -6 $&$ 2 $\\
                   \hline
                   $9 $  &$  $&$  $&$  $&$  $&$  $&$  $&$  $&$  $&$ 1 $\\
                    \hline
 \end{tabular}
 \caption{The plum-blossom product table}
\label{tab:jiandaoji}
\end{table}

Since the plum-blossom product satisfies the commutativity, i.e. $a\clubsuit b=b\clubsuit a$, if all blanks in above table are filled in completely, it can be seen that the whole table is completely symmetric with respect to its main diagonal. It is due to this symmetry that we have left half of our table blank. The advantage of doing so, of course, is that it makes the table more cleaner. It is seen that the third of $45$ plum-blossom products is negative, which ensures that there is a good chance that the pluses and minuses will balance out when calculating the sum of the cross scissor products and the carry from the next position.

\begin{thm}\label{thm-clubsuit}
Given two digits $a\leq b$, suppose that $a\times b=(J(a\clubsuit b),a\clubsuit b)$. Then
\begin{equation*}
J(a\clubsuit b)=
\begin{cases}
a & \textrm{if}\; a=1 \;\textrm{or} \; b=9, \textrm{and} \; b-a\ge 3;\\
a & \textrm{if}\; b-a\ge 5;\\
a-2 &\textrm{if}\; 3\le a\le b\le 7 \;\textrm{and}\; b-a\le 1;\\
a-1 & \;\textrm{otherwise}.
\end{cases}
\end{equation*}
\end{thm}

This theorem incorporating into the basic formula of the rapid multiplication provides a novel multiplication method, which we called the plum-blossom product method. A great amount of examples shows that in most cases the absolute value of the sum of the cross scissor products and the carry from the next position would not be too large, which provides a great convenience for mental multiplication. Let's  see two examples.

\begin{exmp}\rm
Compute $386\times 47. $
\end{exmp}
\par
\solution \hskip 8pt
First by the basic formula of the rapid multiplication, we have $$\begin{bmatrix} 386\\47\end{bmatrix}=(\begin{bmatrix}3\\4 \end{bmatrix},\begin{bmatrix}3 & 8 \\ 4 &7 \end{bmatrix},\begin{bmatrix} 8 & 6 \\ 4 &7 \end{bmatrix}, \begin{bmatrix}6\\7 \end{bmatrix})$$
Then by the plum-blossom product method, we obtain
\begin{align*}
  & 386\times 47\\
=&(3\times 4+3+4-2,3\clubsuit 7+8\clubsuit 4+4+7-2,8\clubsuit 7+6\clubsuit 4+6-1-1,6\clubsuit 7)\\
=&(17,12,-6,2)\\
=&18142.
\end{align*}
\endsolution

\begin{exmp}\rm
Compute $456\times 789. $
\end{exmp}
\par
\solution \hskip 8pt
First by the basic formula of the rapid multiplication, we have
$$\begin{bmatrix} 456\\789\end{bmatrix}
=(\begin{bmatrix}4\\7 \end{bmatrix},\begin{bmatrix}4 & 5 \\ 7 & 8 \end{bmatrix},\begin{bmatrix}4 & 5 & 6\\ 7 & 8 & 9\end{bmatrix},
\begin{bmatrix} 5 & 6 \\ 8 & 9 \end{bmatrix}, \begin{bmatrix}8 \\ 9 \end{bmatrix})$$
Then by the plum-blossom product method, we obtain
\begin{align*}
  & 456\times 789\\
=&(4\times 7+4+5-2,4\clubsuit 8+5\clubsuit 7+4+5+6-3,4\clubsuit 9+5\clubsuit 8+6\clubsuit 7+5+6-2+1,\\
 & 5\clubsuit 9+6\clubsuit 8,6\times 9)\\
=&(35,9,8,-2,4)\\
=&359\,784.
\end{align*}
\endsolution

The plum-blossom product method of multiplication will be applied to the vertical multiplication formula in the next section.

\section{The application of the vertical multiplication formula to mental arithmetic}

We can use the vertical multiplication formula to oral multiplication of two multi-digit numbers, where we could use plum-blossom product method to simplify the computation of vertical multiplication as well as the tare.
We present some examples as follows.

\begin{exmp}\rm
Compute $67\times 89. $
\end{exmp}
\par
\solution \hskip 8pt
The difference product is $\begin{bmatrix}6-7\\8-9 \end{bmatrix}=1$. Then the tare is $10$. The vertical product is $(\begin{bmatrix}6\\8 \end{bmatrix}, \begin{bmatrix}7\\9 \end{bmatrix})=(48+6,7\clubsuit 9)=543=6\overline{6}3$. Thus
$6\overline{6}3\times 11=(6,6+\overline{6},\overline{6}+3,3)=60\overline{3}3=5973$; subtracting the tare  $10$, this yields the answer $5963$.
\endsolution

\begin{exmp}\rm
Compute $677\times 338. $
\end{exmp}
\par
\solution \hskip 8pt
The difference product sums are successively  $K_{12}=\begin{bmatrix}6-7\\3-3 \end{bmatrix}=0$,
$K_{123}=\begin{bmatrix}6-7\\3-8 \end{bmatrix}=5$,
$K_{23}=\begin{bmatrix}7-7\\3-8 \end{bmatrix}=0$. Thus the tare is
$K=(0,5,0,0)=(0,0,0,5,0,0).$
The vertical product is
$$(\begin{bmatrix}6\\3 \end{bmatrix}, \begin{bmatrix}7\\3 \end{bmatrix}, \begin{bmatrix}7\\8 \end{bmatrix})\\
=(6\times 3+2,7\clubsuit 3+5,6)=(2, 0, 6, 6)=(2,1,\overline{4},6),$$
multiplied by $111$ which deduces
$$(2,2+1,2+1+\overline{4}, 1+\overline{4}+6, \overline{4}+6, 6)=(2,3,\overline{1},3,2,6);$$
subtracting the tare $K=(0,0,0,5,0,0)$, this yields that
$$(2,3,\overline{1},3,2,6)-(0,0,0,5,0,0)=(2,3,\overline{1},\overline{2},2,6)=(2,2,8,8,2,6).$$  Therefore,
$677\times 338=228\,826$.
\endsolution

\begin{exmp}\rm
Compute $6789\times 6789. $
\end{exmp}
\par
\solution \hskip 8pt
The difference product sums are successively  $K_{12}=\begin{bmatrix}6-7\\6-7 \end{bmatrix}=1$,
$K_{123}=\begin{bmatrix}6-8\\6-8 \end{bmatrix}=4$,
$K_{1234}=\begin{bmatrix}6-9\\6-9 \end{bmatrix}+\begin{bmatrix}7-8\\7-8 \end{bmatrix}=10$,
$K_{234}=\begin{bmatrix}7-9\\7-9 \end{bmatrix}=4$,
$K_{34}=\begin{bmatrix}8-9\\8-9 \end{bmatrix}=1$. Thus the tare is
$$K=(0,1,4,10,4,1,0)=(0,1,5,0,4,1,0).$$
The vertical product is
$$(\begin{bmatrix}6\\6 \end{bmatrix}, \begin{bmatrix}7\\7 \end{bmatrix}, \begin{bmatrix}8\\8 \end{bmatrix}, \begin{bmatrix}9\\9 \end{bmatrix})\\
=(36+5,7\clubsuit 7+7,8\clubsuit 8+8,9\clubsuit 9)=(41,6,2,1)=(4,2,\overline{4},2,1),$$
multiplied by $1111$ which deduces that
$$(4,4+2,4+2+\overline{4},4+2+\overline{4}+2, 2+\overline{4}+2+1, \overline{4}+2+1, 2+1,1)=(4,6,2,4,1,\overline{1},3,1);$$
subtracting the tare $K=(0,0,1,5,0,4,1,0)$, this yields that
$$(4,6,2,4,1,\overline{1},3,1)-(0,0,1,5,0,4,1,0)=(4,6,1,\overline{1},1,\overline{5},2,1)=(4,6,0,9,0,5,2,1).$$  Therefore,
$6789\times 6789=46\,090\,521$.
\endsolution

\begin{exmp}\rm
Compute $4657\times 86. $
\end{exmp}
\par
\solution \hskip 8pt
We regard $86$ as 4-digit number $0086$ so that $4657\times 86=4657\times 0086$.\\
The difference product sums are successively  $K_{12}=\begin{bmatrix}4-6\\0-0 \end{bmatrix}=0$,
$K_{123}=\begin{bmatrix}4-5\\0-8 \end{bmatrix}=8$,
$K_{1234}=\begin{bmatrix}4-7\\0-6 \end{bmatrix}+\begin{bmatrix}6-5\\0-8 \end{bmatrix}=10$,
$K_{234}=\begin{bmatrix}6-7\\0-6 \end{bmatrix}=6$,
$K_{34}=\begin{bmatrix}5-7\\8-6 \end{bmatrix}=\overline{4}$. Thus the tare is
$$K=(0,8,10,6,\overline{4},0)=(0,9,0,5,6,0)=(1,\overline{1},1,\overline{5},6,0).$$
The vertical product is
$$(\begin{bmatrix}4\\0 \end{bmatrix}, \begin{bmatrix}6\\0 \end{bmatrix}, \begin{bmatrix}5\\8 \end{bmatrix}, \begin{bmatrix}7\\6 \end{bmatrix})\\
=(0,0+4,5\clubsuit 8+4,7\clubsuit 6)=(4,4,2),$$
multiplied by $1111$ which deduces that
$$=(4,8,t,t,6,2)=(4,9,1,0,6,2);$$
subtracting the tare $K=(0,0,1,5,0,4,1,0)$, this yields that
$$(4,9,1,0,6,2)-(1,\overline{1},1,\overline{5},6,0)=(3,t,0,5,0,2)=(4,0,0,5,0,2).$$ Therefore,
$4657\times 86=400\,502$.
\endsolution

\begin{exmp}\rm
Compute $268\times 47. $
\end{exmp}
\par
\solution \hskip 8pt
We regard $268$ as segmented number $(24,28)$ so that $268\times 47=(24,28)\times (4,7)$.\\
The difference product is  $K_{12}=\begin{bmatrix}24-28\\4-7 \end{bmatrix}=12$. Thus the tare is $K=120.$
The vertical product is
$$(\begin{bmatrix}24\\4 \end{bmatrix}, \begin{bmatrix}28\\7 \end{bmatrix})=(96,196)=(10\overline{4},20\overline{4})=(10\overline{4}+20,\overline{4})=12\overline{4}\overline{4}=1156,$$
multiplied by $11$ which deduces that：$$(1,2,6,11,6);$$
subtracting the tare $K=120$, this yields that
$$(1,2,5,9,6).$$  Therefore, $268\times 47=12\,596$.
\endsolution

\begin{exmp}\rm
Compute $29\times 86. $
\end{exmp}
\par
\solution \hskip 8pt
We rewrite $29$ as $3\overline{1}$ so that $29\times 86=3\overline{1}\times 86$.
The difference product sum is  $K_{12}=\begin{bmatrix}3-\overline{1}\\8-6 \end{bmatrix}=8$. Thus the tare is $K=80.$
The vertical product is
$$(\begin{bmatrix}3\\8 \end{bmatrix}, \begin{bmatrix}\overline{1}\\6 \end{bmatrix})=(24,\overline{6})=234,$$
multiplied by $11$ which deduces $2574$;
subtracting the tare $K=80$, this yields that
$$(2,5,7,4)-(0,0,8,0)=(2,5,\overline{1},4)=(2,4,9,4).$$  Therefore, $29\times 86=2\,494$.
\endsolution

\begin{exmp}\rm
Compute $6\,162\times 8\,384. $
\end{exmp}
\par
\solution \hskip 8pt
Using the segment length $2$, we have $6\,162\times 8\,384=(61,62)\times (83,84)$.\\
The difference product sum is  $K_{12}=\begin{bmatrix}61-62\\83-84 \end{bmatrix}=1$. Thus the tare is $K=(1,0).$
Calculating by the plum-blossom product method, we get
\begin{gather*}
61\times 83=(6\times 8+2+1,1\clubsuit 8+6\clubsuit 3,1\times 3)=(51,-4,3)=(51,\overline{4}3),\\
62\times 84=(6\times 8+4+1,2\clubsuit 8+6\clubsuit 4,2\times 4)=(53,-10,8)=(52,08).
\end{gather*}
Thus the vertical product is
$$(\begin{bmatrix}61\\83 \end{bmatrix}, \begin{bmatrix}62\\84 \end{bmatrix})=(51,\overline{4}3+52,08)=(51,15,08),$$
multiplied by $(1,1)=(01,01)$ which deduces that $$(51,51+15,15+08,08)=(51,66,23,08);$$
subtracting the tare $K=(1,0)=(01,00)$, this yields that
$$(51,66,23,08)-(01,00)=(51,66,22,08).$$  Therefore, $6\,162\times 8\,384=51\,662\,208. $
\endsolution

\section{The application of the vertical multiplication formula to computer science}

Another important application of the vertical multiplication formula is in computer science. By this formula, we can give a new multiplication algorithm of long integers. This algorithm is similar to the Karatsuba algorithm \cite{Karatsuba,Rafferty,Sadofsky,Harvey-Hoeven}.

Given two $n$-digit binary numbers $a$ and $b$,  let $C=(C_1,C_2,\ldots,C_n)$ be the vertical product of $a$ and $b$. Then \[C\times (\underbrace{1,1,\ldots,1}_n)=C\times (1,\underbrace{0,\ldots,0}_{n-1},-1)=C\times (1,-1)_n=(C,-C)_n, \]
where $()_n$ denotes the segmented number with each block length $n$. Let $K$ be the tare of $a\times b$. Then the vertical multiplication formula in binary system becomes
\[a\times b=(C,-C)_n-K.\]
For any two $n$-bit numbers, we may directly use this formula to compute the product $a\times b$, see the following Algorithm 1.

\begin{table}[!htb]
\centering
 \begin{tabular} %%|p{1cm}<{\centering}|
 {l}
                    \hline
                  {\bf Algorithm 1:} Vertical Multiplication\\
                  \hline
                  \hspace{1em}{\bf Input:} $n$-bit integers $a$ and $b$\\
                  \hspace{1em}{\bf Output:} $z=a\times b$  \\
                  \hspace{1.5em}1: {\bf for} $i$ in $0$ to $2n-4$ {\bf do}\\
                  \hspace{1.5em}2: \hspace{1em}{\bf if} $i<n-1$ {\bf then}\\
                  \hspace{1.5em}2: \hspace{2em}$K_i=\sum_{j=0}^{[\frac{i}{2}]+1}(a_j-a_{i+1-j})(b_j-b_{i+1-j})$\\
                  \hspace{1.5em}2: \hspace{1em}{\bf else}\\
                  \hspace{1.5em}2: \hspace{2em}$K_i=\sum_{j=[\frac{i+1}{2}]+1}^{n-1}(a_j-a_{i+1-j})(b_j-b_{i+1-j})$\\
                  \hspace{1.5em}2: \hspace{1em}{\bf end if}\\
                  \hspace{1.5em}1: {\bf end for}\\
                  \hspace{1.5em}1: {\bf for} $i$ in $0$ to $n-1$ {\bf do}\\
                  \hspace{1.5em}2: \hspace{1em}$C_i=a_i\times b_i$\\
                  \hspace{1.5em}1: {\bf end for}\\
                  \hspace{1.5em}2: \hspace{1em}$z_i=\sum_{i=0}^{n-1}(C_i\ll2^{n+i})-\sum_{i=0}^{n-1}(C_i\ll2^i)-\sum_{i=0}^{n-1}(K_i\ll2^{i+1})$\\
                  \hspace{1em}{\bf return} z\\
                  \hline
 \end{tabular}
 %\caption{The plum-blossom product table}
%\label{tab:jiandaoji}
\end{table}

As mentioned before, the vertical multiplication formula are suitable for the segmented numbers. For example, Algorithm 2 and 3 are suitable for any two $2m$-bit and $3m$-bit integers respectively with each block length $m$. Note that for any integer, some zeros may be supplemented before it if necessary.

\begin{table}[!htb]
\centering
 \begin{tabular} %%|p{1cm}<{\centering}|
 {l}
                    \hline
                  {\bf Algorithm 2:} Vertical Multiplication with $2m$-bit Integers\\
                  \hline
                  \hspace{1em}{\bf Input:} $2m$-bit integers $a$ and $b$, where $a=(a_1,a_0)$ and $b=(b_1,b_0)$ \\
                  \hspace{4.5em}with $a_i, b_j$ being $m$-bit integers\\
                  \hspace{1em}{\bf Output:} $z=a\times b$  \\
                  \hspace{1.5em}1: $z_0=a_0\times b_0$\\
                  \hspace{1.5em}1: $z_2=a_1\times b_1$\\
                  \hspace{1.5em}1: $z_1=z_0+z_1-(a_1-a_0)\times (b_1-b_0)$\\
                  \hspace{1.5em}1: $z=(z_2,z_1,z_0)=\sum_{i=0}^2(z_i\ll2^{in})$\\
                  \hspace{1em}{\bf return} z\\
                  \hline
 \end{tabular}
 %\caption{The plum-blossom product table}
%\label{tab:jiandaoji}
\end{table}
Obviously, Algorithm 2 is  similar with the algorithm for Karattsuba multiplication in \cite{Rafferty}.
\begin{table}[!htb]
\centering
 \begin{tabular} %%|p{1cm}<{\centering}|
 {l}
                    \hline
                  {\bf Algorithm 3:} Vertical Multiplication with $3m$-bit Integers\\
                  \hline
                  \hspace{1em}{\bf Input:} $3m$-bit integers $a$ and $b$, where $a=(a_2,a_1,a_0)$ and $b=(b_2,b_1,b_0)$ \\
                  \hspace{4.5em}with $a_i, b_j$ being $m$-bit integers\\
                  \hspace{1em}{\bf Output:} $z=a\times b$  \\
                  \hspace{1.5em}1: $K_0=(a_1-a_0)\times (b_1-b_0)$\\
                  \hspace{1.5em}1: $K_1=(a_2-a_0)\times (b_2-b_0)$\\
                  \hspace{1.5em}1: $K_2=(a_2-a_1)\times (b_2-b_1)$\\
                  \hspace{1.5em}1: {\bf end for}\\
                  \hspace{1.5em}1: {\bf for} $i$ in $0$ to $n-1$ {\bf do}\\
                  \hspace{1.5em}2: \hspace{1em}$C_i=a_i\times b_i$\\
                  \hspace{1.5em}1: {\bf end for}\\

                  \hspace{1.5em}1: $z_0=C_0$\\
                  \hspace{1.5em}1: $z_1=C_0+C_1-K_0$\\
                  \hspace{1.5em}1: $z_2=C_0+C_1+C_2-K_1$\\
                  \hspace{1.5em}1: $z_3=C_1+C_2-K_2$\\
                  \hspace{1.5em}1: $z_4=C_2$\\
                  \hspace{1.5em}1: $z=(z_4,z_3,z_2,z_1,z_0)=\sum_{i=0}^4(z_i\ll2^{in})$\\

                  \hspace{1em}{\bf return} z\\
                  \hline
 \end{tabular}
 %\caption{The plum-blossom product table}
%\label{tab:jiandaoji}
\end{table}

Now we give a real example of calculating the product of two binary numbers by the vertical multiplication formula.
\begin{exmp}\rm
$111101\times 101011=?$
\end{exmp}
\par
\solution \hskip 8pt
{\bf Method 1. }
 First, we calculate the difference products as follows:
\begin{gather*}
K_{12}=K_{13}=K_{14}=K_{15}=K_{16}=0, \\
K_{23}=K_{24}=K_{26}=0, K_{25}=-1, \\
K_{34}=K_{35}=K_{36}=0, K_{45}=-1, K_{46}=0, K_{56}=0.
\end{gather*}
Then the symmetric difference product sums are in turn
\begin{gather*}
K_{12}=0, K_{123}=K_{13}=0, K_{1234}=K_{14}+K_{23}=0,\\
K_{12345}=K_{15}+K_{24}=0, K_{123456}=K_{16}+K_{25}+K_{34}=-1, \\
K_{23456}=K_{26}+K_{35}=0, K_{3456}=K_{36}+K_{45}=-1, K_{456}=K_{46}=0, K_{56}=0.
\end{gather*}
Thus the tare is
\begin{align*} K &=(K_{12},K_{123},K_{1234},K_{12345}, K_{123456}, K_{23456}, K_{3456}, K_{456}, K_{56}, 0)\\
                &=(0,0,0,0,-1,0,-1,0,0,0)=(-1,0,-1,0,0,0).
\end{align*}
At last, the vertical product is as follows:
\begin{equation*}
C=(C_1,C_2,C_3,C_4,C_5,C_6)=(1,0,1,0,0,1).
\end{equation*}
By the vertical multiplication formula, we have that
\begin{align*}
  & 111101\times 101011 \\
 =& (1,0,1,0,0,1,-1,0,-1,0,0,-1)-(0,0,0,0,-1,0,-1,0,0,0)\\
 =& (1,0,1,0,0,1,-1,0,-1,0,0,-1)+(1,0,1,0,0,0)\\
 =& (1,0,1,0,0,1,0,0,0,0,0,-1)\\
 =& (1,0,1,0,0,0,1,1,1,1,1,1)\\
 =& 101000111111.
\end{align*}

{\bf Method 2. }
we use the numbers with segment length $2$. Then
\[111101\times 101011=(11,11,01)\times (10,10,11).\]
The difference products are as follows:
\begin{gather*}
K_{12}=0, K_{13}=K_{23}=-10.
\end{gather*}
Then the symmetric difference product sums are in turn
\begin{gather*}
K_{12}=0, K_{123}=K_{13}=-10, K_{23}=-10.
\end{gather*}
Thus the tare is
\begin{align*} K &=(K_{12},K_{123},K_{23}, 0)\\
                &=(00,-10,-10,00)=-(10,10,00).
\end{align*}
At last, the vertical product is as follows:
\begin{align*}
C &=(C_1,C_2,C_3)\\
&=(11\times 10,11\times 10,01\times 11)\\
&=(110,110,11)=(1,11,10,11).
\end{align*}
By the vertical multiplication formula, we have that
\begin{align*}
  & 111101\times 101011 =(11,11,01)\times (10,10,11) =C\times (01,01,01)-K\\
 =& (1,11,10,11)\times (01,01,01)+(10,10,00)\\
 =& (1,1+11,1+11+10,11+10+11,10+11,11)+(10,10,00)\\
 =& (1,1+11,1+11+10,11+10+11+10,10+11+10,11)\\
 =& (10,10,00,11,11,11)\\
 =& 101000111111.
\end{align*}

{\bf Method 3. }
we use the numbers with segment length $3$. Then
\[111101\times 101011=(111,101)\times (101,011).\]
The difference product is as follows:
\begin{gather*}
K_{12}=(111-101)\times (101-011)=10\times 10=100.
\end{gather*}
Then  the tare is
\begin{align*} K &=(K_{12}, 0)\\
                &=(100,000).
\end{align*}
At last, the vertical multiplication is as follows:
\begin{align*}
C &=(C_1,C_2)\\
&=(111\times 101,101\times 011)\\
&=(100011,1111)=(100,100,111).
\end{align*}
By the vertical multiplication formula, we have that
\begin{align*}
  & 111101\times 101011 =(111,101)\times (101,011) =C\times (001,001)-K\\
 =& (100,100,111)\times (001,001)-(100,000)\\
 =& (100,100+100,100+111,111)-(100,000)\\
 =& (100,100+100,111,111)\\
 =& (101,000,111,111)\\
 =& 101000111111.
\end{align*}
\endsolution

Let's turn to the analysis of the complexity of our algorithm.

If an addition or a multiplication of tow digits is called one step of an algorithm, then the time complexity of the algorithm may be calculated by the whole steps.
Let $a,b$ be $n$-bit binary numbers with $n=k^m$. Then $m=\log_k n$. We divide $a,b$ into $k$'s blocks with each block length $b=k^{m-1}$. We denote $b=(b_1,b_2,\ldots,b_k)$ and $b=(b_1,b_2,\ldots,b_k)$. Then by $T(n)$, we denote the complexity of calculating the product $a\times b$ by the vertical multiplication. The amount of the difference products of $a_i$ and $b_j$ is $C_k^2=k(k-1)/2$, which plus the amount of the vertical products $C_i=a_i\times b_i$ is equal to $k(k+1)/1$. Each difference product contains $1$ subtraction, and $C\times (\underbrace{1,\ldots,1}_k)$ contains $2bk=2n$ additions of digits. The difference product sum contains $\frac{k^2}{4}$ additions of $b$-digits numbers, that is, it contains $\frac{k^2\cdot 2b}{4}=\frac{kn}{2}$ additions of digits.
Therefore we have the following recursive formula:
\begin{equation}\label{eq.1}
T(n)=\frac{k(k+1)}{2}\cdot T(\frac{n}{k})+\frac{kn}{2}+2n.\end{equation}
When $k=2$, the above recursive formula becomes
\[ T(n)=3\cdot T(\frac{n}{2})+3n,\]
which deduces that
\begin{equation} T(n)=7n^{\log_2 3}\thicksim n^{1.58};\end{equation}
When $k=3$, the recursive formula (\ref{eq.1}) becomes
\[ T(n)=6\cdot T(\frac{n}{3})+\frac{7n}{2},\]
which deduces that
\begin{equation} T(n)=\frac{9}{2}\cdot n^{\log_3 6}\thicksim n^{1.63}.\end{equation}
Note that when $k=2$, the complexity of the vertical multiplication is equals to that of Karatsuba multiplication.

%%%%%%%%%%%

\section{Conclusion}

As we have seen above, our vertical multiplication formula can be used not only to design the multiplier but also to mental multiplication.
Compared with the Vedic multiplications, our vertical multiplication formula provides a more general and much systematic mental multiplication method.
Our algorithm is a generalization of that of Karatsuba and Offman, but it is  superior to the latter because the latter is not suitable for oral calculation of multiplication. Two algorithms are both the vertical multiplication method, but in our algorithm the horizontal subtraction in stead of the horizontal addition used in Kartsuba algorithm makes the involved numbers are much smaller and their algebraic sum are much smaller too because the positive and negative numbers probably cancel each other out. If incorporating into other methods such as FFT and the modular method, our vertical multiplication formula will allow more practical and efficient multipliers to be designed.

%%\section*{Acknowledgement} The author thanks the editor for his helpful suggestions!

\end{sloppypar}  %%% %解决行溢出的问题 solvee overfull\hbox

\end{document}